\documentclass[reqno, hyp]{nyj}
\title{Non-ergodic actions, cocycles and superrigidity}

\author{David Fisher}
 \address{Department of Mathematics and Computer Science,
 Lehman College -- CUNY,
 250 Bedford Park Boulevard W.,
 Bronx, NY 10468}
 \email{dfisher@lehman.cuny.edu,
  http://comet.lehman.cuny.edu/fisher/}

\author{Dave Witte Morris}
 \address{Department of Mathematics and Computer Science,
 University of Lethbridge,
 4401 University Drive,
 Lethbridge, AB, T1K~3M4, Canada}
 \email{Dave.Morris@uleth.ca,
 http://people.uleth.ca/$\sim$dave.morris/}

\author{Kevin Whyte}
 \address{Department of Mathematics,
 University of Illinois at Chicago,
 851 S.~Morgan Street,
 Chicago, IL 60607--7045}
 \email{kwhyte@math.uic.edu,
 http://www2.math.uic.edu/$\sim$kwhyte/}

\thanks{This research was partially supported by the National Science
Foundation (grants DMS--0226121, DMS--0100438, and DMS--0204576).
The first author was also partially supported by a PSC-CUNY
grant.}



\newcommand{\symmdiff}{\mathbin{\Delta}}
\DeclareMathOperator{\Id}{Id}
\newcommand{\real}{\mathbb{R}}
\renewcommand{\natural}{\mathbb{N}}
\newcommand{\func}{\mathcal{F}}
\newcommand{\F}{\mathbf{F}}
\newcommand{\erg}{\Omega}
\newcommand{\cobound}{\delta}
\DeclareMathOperator{\Aut}{Aut}
\newcommand{\A}{\Aut_{[\mu]}(S)}
\newcommand{\B}{\mathcal{B}}
\DeclareMathOperator{\cpct}{Cpct}
\DeclareMathOperator{\dist}{dist}
\DeclareMathOperator{\diam}{diam}
\newcommand{\triv}{\mathord{\bf 1}}

\DeclareMathOperator{\Const}{Const}
\DeclareMathOperator{\EssRg}{EssRg}
\DeclareMathOperator{\Stab}{Stab}
\newcommand{\field}{\mathbb{F}}
\newcommand{\proj}[1]{\mathbb{P}(#1)}
\newcommand{\Unitary}{\mathord{\mathbf{U}\bigl( L^2(S) \bigr)}}
\DeclareMathOperator{\Prob}{Prob}

  \newcommand{\G}{\Gamma}
  \newcommand{\g}{\gamma}
  \newcommand{\Ga}{\mathbb G}

\newcommand{\bigset}[2]{\left\{\, #1
 \mathrel{\left| \vphantom {\left\{ #1 \mid #2 \right\} }
 \right.} #2 \,\right\} }

\renewcommand{\see}[1]{{\upshape(}see~\ref{#1}{\upshape)}}

\newcommand{\cf}[1]{(cf.~\ref{#1})}
\newcommand{\pref}[1]{{\upshape(}\ref{#1}{\upshape)}}
\newcommand{\fullref}[2]{\ref{#1}\pref{#1-#2}}
\newcommand{\fullsee}[2]{{\upshape(}see~\ref{#1}\pref{#1-#2}{\upshape)}}

\renewcommand{\thesubsection}{\thesection\Alph{subsection}}

\numberwithin{equation}{section}

 \newtheorem{thm}[equation]{Theorem}
 \newtheorem{lem}[equation]{Lemma}
 \newtheorem{prop}[equation]{Proposition}
 \newtheorem{cor}[equation]{Corollary}

\theoremstyle{definition}
 \newtheorem{defn}[equation]{Definition}
 \newtheorem{rem}[equation]{Remark}
 \newtheorem{notn}[equation]{Notation}

 \newtheorem{ack}{Acknowledgments\ignorespaces}

 \newcounter{steppf}

\newenvironment{steppf}[1][\unskip]{\refstepcounter{steppf}
 \em
 \medskip \noindent Step \thestep\ #1.\ }{\unskip\upshape}
 \newcommand{\thestep}{\arabic{steppf}}

\begin{document}

\begin{abstract}
This paper proves various results concerning non-ergodic actions
of locally compact groups and particularly Borel cocycles defined
over such actions. The general philosophy is to reduce the study
of the cocycle to the study of its restriction to each ergodic
component of the action, while being careful to show that all
objects arising in the analysis depend measurably on the ergodic
component.  This allows us to prove a version of the superrigidity
theorems for cocycles defined over non-ergodic actions.
\end{abstract}

\maketitle

\tableofcontents

\section{Introduction}

It is often the case that one has extensive information about each
ergodic component of an action, and would like to piece this local
information together (measurably) in order to obtain a global
conclusion about the entire action. This note addresses a number of
problems of this type, mostly dealing with cocycles. For example:
 \begin{enumerate}
 \item If $\alpha$ and $\beta$ are Borel cocycles, and the restriction
of~$\alpha$ to almost every ergodic component is cohomologous to the
restriction of~$\beta$, then $\alpha$ is cohomologous to~$\beta$
\see{CohoErg}.
 \item  If $\alpha$ is a Borel cocycle, and the restriction of~$\alpha$
to almost every ergodic component is cohomologous to a homomorphism
cocycle, then $\alpha$ is cohomologous to a homomorphism cocycle
\see{CohoHomo}.
 \item If almost every ergodic component of a $G$-action has a certain
standard Borel $G$-space~$X$ as a measurable quotient, then $X$~is a
measurable quotient of the entire action \see{QuotErgComps}.
 \end{enumerate}
 (We consider only Borel actions of second countable, locally compact
groups on standard Borel spaces with a quasi-invariant probability
measure.)

We also prove a superrigidity theorem for cocycles that applies to
non-ergodic actions \see{theorem:Gsuperrigidityne}.  In fact our
work was motivated by the discovery, during the writing of
\cite{FisherMargulis}, that no proof of any version of
superrigidity for cocycles concerning cocycles over non-ergodic
actions, exists in the literature.  For many applications to
non-ergodic actions, including those in \cite{FisherMargulis},
other information concerning the action and cocycle allows one to,
with some additional work make do with superrigidity theorems for
cocycles which are defined over ergodic actions.  However, the
results in section \ref{section:superrigid} allow some
simplification of the arguments in section $5$ of
\cite{FisherMargulis}, and should have other applications as well.

\begin{ack}
 Most of this research was accomplished during a visit to the Institute
for Mathematical Research (FIM) of the Swiss Federal Institute of
Technology (ETH) in Zurich. We are pleased to thank FIM for its
hospitality, and for the financial support that made the visit possible.
We would also like to thank Scot Adams, A.~Kechris, and G.~A.~Margulis
for helpful comments during the course of the investigation.
 \end{ack}

\section{Some lemmas on measurability} \label{LemmaSect}

This section records basic definitions and notation, and also proves that
various natural constructions of sets, functions, and actions yield
results that are measurable. Some of the conclusions are known (or even
well known), but others may be of independent interest.

\subsection{Properties of standard Borel spaces}

We assume the basic theory of Polish spaces, standard Borel spaces and
analytic Borel spaces, which can be found in a number of textbooks, such
as \cite[Chap.~3]{Arveson}. We recall the definition of a standard space
and an analytic set.

\begin{defn} \
 \begin{itemize}
 \item A topological space is \emph{Polish} if it is homeomorphic to
a complete, separable metric space.
 \item A Borel space is \emph{standard} if it is Borel isomorphic to
a Polish topological space.
 \item The pair $(S,\mu)$ is a \emph{standard Borel probability space} if
$S$~is a standard Borel space and $\mu$ is a probability measure on~$S$.
 \item A subset~$A$ of a standard Borel space $S$ is \emph{analytic} if
there exist
 \begin{itemize}
 \item a standard Borel space~$X$
 and
 \item a Borel map $\psi \colon X \to S$,
 \end{itemize}
 such that $A = \psi(X)$ is the image of~$\psi$.
 \end{itemize}
 \end{defn}

\begin{rem}[{\cite[Thm.~3.2.4, p.~67]{Arveson}}] \label{Anal->Meas}
 Any analytic subset~$A$ of any standard Borel space~$S$ is
\emph{absolutely measurable}; that is, for any probability measure~$\mu$
on~$S$, there exist Borel subsets $B_1$ and~$B_2$ of~$S$ with $B_1
\subset A \subset B_2$ and $\mu(B_2 \smallsetminus B_1) = 0$.
 \end{rem}

\begin{thm}[{von Neumann Selection Theorem, cf.\  \cite[Thm.~3.4.3,
p.~77]{Arveson}}] \label{vonNeumann}
 Let
 \begin{itemize}
 \item $(\erg,\nu)$ be a standard Borel probability space,
 \item $L$ be a standard Borel space,
 \item $\func$ be an analytic subset of $\erg \times L$,
 and
 \item $\erg_{\func}$ be the projection of~$\func$ to~$\erg$.
 \end{itemize}
 Then there are
 \begin{itemize}
 \item a conull, Borel subset~$\erg_0$ of~$\erg_{\func}$,
 and
 \item a Borel function $\Phi \colon \erg_0 \to L$,
 \end{itemize}
 such that $\bigl( \omega, \Phi(\omega) \bigr) \in \func$, for all
$\omega \in \erg_0$.
 \end{thm}

\begin{notn}
 Suppose $(S,\mu)$ is a standard Borel probability space, and
$(X,d)$~is a separable metric space.
 \begin{enumerate}
 \item We use $\F(S,X)$ to denote the space of
measurable functions $f \colon S \to X$, where two functions are
identified if they are equal almost everywhere. If $X$~is complete, then
$\F(S,X)$ is a Polish space, under the topology of convergence in
measure (cf.\  \cite[\S4.4]{WheedenZygmund}). A metric can be given by
 $$ d_{\F}(f,g) = \min \bigset{ \vphantom{\Bigl|} \epsilon \ge 0}{
 \mu \bigl\{\, s \in S \mid d \bigl( f(s), g(s) \bigr) > \epsilon
\,\bigr\} \le \epsilon} $$

 \item We use $\B(S)$ to denote the Boolean algebra of measurable subsets
of~$S$, where two subsets are identified if their symmetric difference
has measure~$0$. It is well known that this is a complete separable metric
space, with metric
 $$ d_{\B}(A,B) = \mu(A \symmdiff B) .$$
 \end{enumerate}
 \end{notn}

\begin{rem}
 The $\sigma$-algebra of Borel subsets of $\F(S,X)$ is generated by the
sets of the form
 $$ \Delta_{S_0,X_0,\epsilon} = \bigset{f }{ \mu \bigl( S_0 \cap
f^{-1}(X_0) \bigr) < \epsilon } ,$$
 where $S_0$ and~$X_0$ are Borel subsets of $S$ and~$X$,
respectively, and $\epsilon > 0$. This implies that if $(X',d')$ is Borel
isomorphic to $(X,d)$, then $\F(S,X')$ is Borel isomorphic to $\F(S,X)$.
 \end{rem}

The following is well known:

\begin{lem}[{cf.\ \cite[Lem.~7.1.3, p.~215]{MargulisBook}}]
\label{FiberFunc}
 Let
 \begin{itemize}
 \item $L$ and~$S$ be standard Borel spaces,
 \item $\mu$ be a probability measure on~$S$,
 \item $X$ be a separable metric space,
 and
 \item $f \colon L \times S \to X$ be Borel.
 \end{itemize}
 Then
 \begin{enumerate}
 \item \label{FiberFunc-fiber}
 for each $\ell \in L$, the function $f_\ell \colon S \to X$, defined by
$f_\ell(s) = f(\ell,s)$, is Borel,
 and
 \item \label{FiberFunc-Borel}
 the induced function $\check f \colon L \to \F(S,X)$, defined by $\check
f(\ell) = f_\ell$, is Borel.
 \end{enumerate}
 \end{lem}

In short, any Borel function from $L \times S$ to~$X$ yields a Borel
function from~$L$ to $\F(S,X)$. The converse is true:

\begin{lem} \label{DefdFiberwise}
 Let
 \begin{itemize}
 \item $L$ and $S$ be standard Borel spaces,
 \item $\mu$ be a probability measure on~$S$,
 \item $(X,d)$ be a separable metric space,
 and
 \item $\phi \colon L \to \F(S,X)$ be a Borel function.
 \end{itemize}
 Then there is a Borel function $\hat\phi \colon L \times S \to X$, such
that, for each $\ell \in L$, we have
 \begin{equation} \label{phi(l,s)=phi(l)(s)}
 \mbox{$\hat\phi(\ell,s) = \phi(\ell)(s)$ for a.e.\ $s \in
S$.}
 \end{equation}
 \end{lem}

\begin{proof}
  For each $n \in \natural$, let $\{D_n^i\}_{i=1}^\infty$ be a partition
of $\F(S,X)$ into countably many (nonempty) Borel sets of diameter
less than $2^{-n}$, and choose $\phi_n^i \in D_n^i$. Define
$\phi_n \colon L \times S \to X$ by
  $$ \mbox{$\phi_n(\ell,s) = \phi_n^i(s)$ if $\phi(\ell) \in D_n^i$.}
$$
  Then each $\phi_n$ is Borel.

By replacing each $\{D_n^i\}_{i=1}^\infty$ by the join of $\{D_1^i\}_{i=1}^\infty,
\ldots, \{D_n^i\}_{i=1}^\infty$ we may assume that $\{D_{n+1}^i\}_{i=1}^\infty$ is
a refinement of $\{D_n^i\}_{i=1}^\infty$.  For each $m,n \in \natural$ and each $\ell \in L$
  $$ \mu \bigl\{\, s \in S \mid d\bigl( \phi_m(\ell,s), \phi_n(\ell,s)
\bigr) > 2^{-\min(m,n)} \,\bigr\} < 2^{-\min(m,n)} .$$
Thus, on each fiber $\{\ell\} \times S$, the sequence $\{\phi_n\}$ not only
converges in measure, but converges quickly. There is no harm in
assuming that $X$~is complete; then $\{\phi_n\}$ converges
pointwise a.e. Let $\hat\phi$ be the pointwise limit of
$\{\phi_n\}$. (Define $\hat\phi$ to be constant on the set where
$\{\phi_n\}$ does not converge.) Then $\hat\phi$ is Borel, and
satisfies \pref{phi(l,s)=phi(l)(s)}.

 \end{proof}

\begin{cor} \label{F(F)}
 Let
 \begin{itemize}
 \item $L$ and $S$ be standard Borel probability spaces,
 and
 \item $(X,d)$ be a separable metric space.
 \end{itemize}
 Then $\F( L \times S, X)$ is naturally homeomorphic to $\F \bigl( L, \F(
S, X) \bigr)$.
 \end{cor}

\begin{proof}
 From the preceding two lemmas, we know that there is a natural bijection
between the two spaces. Fubini's Theorem implies that a sequence
converges in one of the spaces if and only if the corresponding sequence
converges in the other space.
 \end{proof}

\begin{defn} Suppose $(S,\mu)$ is a standard Borel probability space.
 \begin{itemize}
 \item Let $\A$ be the group of all equivalence classes of measure-class-preserving Borel
automorphisms of $(S,\mu)$, where two automorphisms are equivalent if they are
equal almost everywhere.
 \item Let $\Unitary$ be the group of unitary operators on the Hilbert
space $L^2(S)$, with the strong operator topology (that is, $T_n \to T$
if $\|T_n(f) - T(f)\| \to 0$, for every $f \in L^2(S)$, equivalently, the topology on
$\Unitary$ has a subbasis of open sets $\mathcal{U}(f,g,\epsilon)=\{T : \|Tf-g\|<\epsilon\}$).  Note that $\Unitary$ is a Polish space.
 \item There is a well known embedding of $\A$ in $\Unitary$, given by
 $$ T_\phi(f)(s) = D_\phi(s)^{1/2} \, f\bigl( \phi^{-1}(s) \bigr) ,$$
 for $\phi \in \A$ and $f \in L^2(S)$, where $D_\phi$~is the Radon-Nikodym
derivative of~$\phi$. This provides $\A$ with the topology
of a separable metric space, and thereby makes $\A$ into a topological
group.
 \end{itemize}
 \end{defn}

\begin{rem}
 Note that $\A$ is \emph{not} locally compact.
 On the other hand, $\A$ is a closed subset of $\Unitary$ (because it
consists of the operators that map nonnegative functions to nonnegative
functions \cite[\S3]{GoodrichEtAl}), so it is a Polish space.
 \end{rem}

\begin{prop}[{Ramsay, cf.\ \cite[Cor.~3.4]{Ramsay}}]
 If $(S,\mu)$ is a standard Borel probability space, then $\A$ acts
continuously on the Borel algebra $\B(S)$.
 \end{prop}

\begin{proof}
 We wish to show that if $\phi$ is close to $\phi_0$ in $\A$, and $A$ is close
to~$A_0$ in $\B(S)$, then $\mu \bigl( \phi(A) \smallsetminus \phi_0(A_0) \bigr)$
and  $\mu \bigl( \phi_0(A_0) \smallsetminus \phi(A) \bigr)$ are close to~$0$.
Thus, letting $\psi$ be $\chi_A$ and $\chi_{A_0}$, it suffices to show that if
$\psi \in L^2(S)$, with $\|\psi\| \le 1$, then
 $\bigl| \int_{\phi(A)} \psi^2 \,d\mu - \int_{\phi_0(A_0)} \psi^2 \,d\mu \bigr| $
 is close to~$0$. To simplify the notation, we replace $\phi$ and~$\phi_0$
by their inverses in the following calculation.

We have
 \begin{align*}
 &\left| \int_{\phi^{-1}(A)} \psi^2 \,d\mu - \int_{\phi_0^{-1}(A_0)}
\psi^2 \,d\mu \right| \\
 &\kern2em =  \left| \int_A T_\phi(\psi)^2 \,d\mu - \int_{A_0}
T_{\phi_0}(\psi)^2 \,d\mu \right| \\
 &\kern2em \le \left| \int_A T_\phi(\psi)^2 \,d\mu - \int_{A}
T_{\phi_0}(\psi)^2 \,d\mu \right| \\
 &\kern4em  + \left| \int_A T_{\phi_0}(\psi)^2 \,d\mu -
\int_{A_0} T_{\phi_0}(\psi)^2 \,d\mu \right| \\
 &\kern2em = \left| \int_A \bigl( T_\phi(\psi) + T_{\phi_0}(\psi)
\bigr)
 \bigl( T_\phi(\psi) - T_{\phi_0}(\psi) \bigr) \,d\mu
\right| \\
 &\kern4em  + \left| \int_A T_{\phi_0}(\psi)^2 \,d\mu -
\int_{A_0} T_{\phi_0}(\psi)^2 \,d\mu \right| \\
 &\kern2em \le 2
 \sqrt{\int_A \bigl( T_\phi(\psi) - T_{\phi_0}(\psi) \bigr)^2 \,d\mu}
 && \begin{matrix}
 \mbox{(H\"older's Inequality} \\
 \mbox{and $\|\psi\| \le
1$)}
 \end{matrix}
 \\
 &\kern4em  + \int_{A \symmdiff A_0} T_{\phi_0}(\psi)^2
\,d\mu
 && \begin{matrix}
 \mbox{(integrand is 0} \\
 \mbox{on $A \cap A_0$)}
 . \end{matrix}
 \end{align*}
 By definition of the topology on $\Unitary$, the first term in the final
expression is small whenever $\phi$ is close to~$\phi_0$ (since $\psi$ is
fixed). Because $T_{\phi_0}(\psi)^2$ is a fixed $L^1$~function, the second
term is small whenever $A \symmdiff A_0$ has sufficiently small measure
\cite[Exer.~6.10(a)]{Rudin-RCAnal}; that is, whenever $A$ is sufficiently
close to~$A_0$ in $\B(S)$.
 Thus,
 $\left| \int_{\phi^{-1}(A)} \psi^2 \,d\mu - \int_{\phi_0^{-1}(A_0)}
\psi^2 \,d\mu \right| $
 is close to~$0$, as desired.
 \end{proof}

The fact that the action on $\B(S)$ is continuous
at the empty set~$\emptyset$ can be restated as follows:

\begin{cor} \label{TendToZero}
 If $\phi_n \to \phi$ in $\A$, and $A_n$ is a sequence of Borel subsets
of~$S$, such that $\mu(A_n) \to 0$, then $\mu \bigl( \phi_n(A_n) \bigr)
\to 0$.
 \end{cor}

In the following result, we assume that $S$ is a separable metric space,
so that we can speak of convergence in measure. This is a very mild
assumption, because any standard Borel space is, by definition, Borel
isomorphic to such a space.

\begin{cor} \label{ConvInMeas}
 Assume that $S$ is a separable metric space.
 If $\phi_n \to \phi$ in $\A$, then $\phi_n \to \phi$ in measure.
 \end{cor}

\begin{proof}
 By considering $\phi_n \circ \phi^{-1}$ (and using the fact that the
action of~$\phi$ on $\B(S)$ is continuous at~$\emptyset$), we may assume
that $\phi = \Id$. Let $\{A^i\}_{i=1}^\infty$ be a partition of~$S$ into
countably many Borel sets, such that $\diam(A^i) < \epsilon$ for each~$i$.
 Let
 $$ \Sigma_n^i = \{\, s \in A^i \mid d \bigl( \phi_n(s), s \bigr) >
\epsilon \,\} .$$ Then, for each fixed~$i$, we have
 $$ \mu \bigl( \phi_n (\Sigma^n_i) \bigr)
 \le \mu \bigl( \phi_n(A^i) \smallsetminus A^i \bigr)
 \le d \bigl( \phi_n(A^i), A^i \bigr)
 \to 0
 \mbox{ as $n \to \infty$}.
 $$
 Thus, $\phi_n (\Sigma_n^i) \to \emptyset$ in $\B(S)$. Because $\A$ acts
continuously on $\B(S)$, this implies that
 $$ \mu(\Sigma_n^i) \to 0 \mbox{ as $n \to \infty$}.$$
 Therefore
 \begin{align*}
 &\lim_{n \to \infty} \mu \{\, s \in S \mid d \bigl( \phi_n(s), s
\bigr) > \epsilon \,\} \\
 &\kern3em = \lim_{n \to \infty} \mu \bigl( \cup_{i=1}^\infty \Sigma_n^i
\bigr) \\
 &\kern3em \le  \inf_{k \in \natural} \lim_{n \to \infty} \mu \Bigl(
 \bigl( \cup_{i=1}^k \Sigma_n^i \bigr) \cup \bigl( S
\smallsetminus (A^1 \cup \cdots \cup A^k) \bigr) \Bigr) \\
 &\kern3em = \lim_{k \to \infty} \mu \bigl( S
\smallsetminus (A^1 \cup \cdots \cup A^k) \bigr) \\
 &\kern3em=  0
 . \end{align*}
 So $\phi_n \to \Id$ in measure.
 \end{proof}

\begin{defn}
 There is a natural action of $\A$ on $\F(S,X)$, defined by $\phi(f) = f
\circ \phi^{-1}$, for $\phi \in \A$ and $f \in \F(S,X)$.
 \end{defn}

\begin{prop} \label{Autmu(S)action}
 If $(S,\mu)$ is a standard Borel probability space, and $(X,d)$ is a
separable metric space, then the natural action of $\A$ on $\F(S,X)$ is
continuous.
 \end{prop}

\begin{proof}
 Suppose $\phi_n \to \phi$ in $\A$, and $f_n \to f$ in $\F(S,X)$.  Note
that
 $$ d (  f_n \phi_n, f \phi )
 \le d (  f_n\phi_n, f\phi_n ) + d (  f \phi_n, f \phi ) .$$
 We have  $d (  f_n\phi_n, f\phi_n ) \to 0$, because
 \begin{align*}
 \mu \bigset{ s }{ d \bigl( f_n\phi_n(s), f \phi_n(s) \bigr)
\bigr)> \epsilon }
 &= \mu \Bigl( \phi_n^{-1} \bigl\{\, s \mid d \bigl( f_n(s), f(s)
\bigr) > \epsilon \,\bigr\} \Bigr) \\
 & \to 0,
 \end{align*}
 using \pref{TendToZero} and the fact that $f_n \to f$ in measure to get
the final limit.

Because $(S,\mu)$ is standard, there is no harm in assuming that $S$ is a
complete, separable metric space. Then, by Lusin's Theorem, there is a
(large) subset~$K$ of~$S$, such that $f$ is uniformly continuous on~$K$.
Let
 $$ A_n = \phi^{-1}(S \smallsetminus K) \cup \phi_n^{-1}(S \smallsetminus
K) .$$
 From \pref{TendToZero}, we see that, by requiring $\mu(K)$ to be
sufficiently large, we may ensure that
 $$ \mu (A_n) < \epsilon .$$
 Choose $\delta > 0$, such that $d \bigl( f(s), f(t) \bigr) < \epsilon$,
for all $s,t \in K$ with $d_S(s,t) < \delta$.
 We have
 \begin{align*}
 \limsup d (  f \phi_n, f \phi )
 &= \limsup \mu \bigl\{\, s \mid d \bigl( f\phi_n(s), f \phi(s) \bigr)
> \epsilon\, \bigr\} \\
 &\le \limsup \mu(A_n) + \lim \mu \bigl\{\, s \mid d \bigl( \phi_n(s),
\phi(s) \bigr) > \delta \,\bigr\}
 \\
 &\le \epsilon + 0
 , \end{align*}
 using \pref{ConvInMeas} to obtain the term``$0$" in the final expression.

 Since $\epsilon > 0$ is arbitrary, we conclude that $d(f_n\phi_n, f\phi)
\to 0$.
 \end{proof}

We will use the following easy observation.

\begin{lem} \label{IntIsCont}
 If $A$ is any Borel subset of any standard Borel probability space
$(S,\mu)$, then the integration functional
 $ I \colon \F(S,\real^{\ge0}) \to \real \cup \{\infty\}$,
 defined by $I(f) = \int_A f\, d\mu$, is Borel.
 \end{lem}

\begin{proof}
 Let $\chi_A$ be the characteristic function of~$A$. It is easy to see that
the map $f \mapsto \chi_A f$ is a continuous function on
$F(S,\real^{\ge0})$, so we may (and will) assume $A = S$.

 For each $n \in \natural$, choose a continuous function $\xi_n \colon
\real^{\ge0} \to [0,n]$, such that
 $$ \xi_n(x) =
 \begin{cases}
 x & \mbox{if $x \le n$}, \\
 0 & \mbox{if $x \ge n+1$}.
 \end{cases}
 $$
 Note, for each $f \in \F(S,\real^{\ge0})$, that the composition $\xi_n
\circ f$ is bounded by~$n$, and $\xi_n \circ f \uparrow f$ pointwise.

Because $\xi_n$ is uniformly continuous (being a continuous function with
compact support), it is easy to see that the map $\F(S,\real^{\ge0}) \to
\F(S,[0,n])$, defined by $f \mapsto \xi_n \circ f$, is continuous.
Furthermore, the restriction of~$I$ to $\F(S,[0,n])$ is continuous. Thus,
the function
 $$ \mbox{$ I_n \colon \F(S,\real^{\ge0}) \to \real^{\ge0}$,
 defined by $I_n(f) = I(\xi_n \circ f)$,}$$
 is continuous.
 The Monotone Convergence Theorem implies that
 $$I(f) = \lim_{n \to \infty} I_n(f) ,$$ so $I$~is a pointwise limit of
continuous functions. Therefore, $I$ is Borel.
 \end{proof}

\subsection{Ergodic decomposition and near actions}
 A proof of the following folklore theorem has been provided by
G.~Greschonig and K.~Scmidt \cite{Greschonig-Schmidt}. The statement here
is slightly stronger than \cite[Thm.~1.1]{Greschonig-Schmidt} (because
Thm.~\fullref{ErgDecompThm}{atom} is more precise than
\cite[Thm.~1.1(3)]{Greschonig-Schmidt}), but it follows immediately from
\cite[Thm.~5.2]{Greschonig-Schmidt}, by letting
 $\erg$ be a $p_*\mu$-conull Borel subset of $\Prob(X)$,
 $X' = p^{-1}(\erg)$,
 $\psi(x) = p(x)$,
 and
 $\xi(\omega) = \omega$.

\begin{notn}
 For any standard Borel space~$X$, we use $\Prob(X)$ to denote the space
of all probability measures on~$X$. It is well known that $\Prob(X)$ is a
standard Borel space, under an appropriate weak$^*$ topology.
 \end{notn}

\begin{thm}[{Ergodic Decomposition \cite{Greschonig-Schmidt}}]
\label{ErgDecompThm}
 Let
 \begin{itemize}
 \item $G$ be a locally compact second countable group,
 \item $(X,\mu)$ be a standard Borel probability space,
 and
 \item $\rho \colon G \times X \to X$ be a Borel action, such that
$\mu$~is quasi-invariant.
 \end{itemize}
 Then there exist:
 \begin{itemize}
 \item a standard Borel probability space $(\erg,\nu)$,
 \item a conull $G$-invariant Borel subset $X'$ of~$X$,
 \item a $G$-invariant Borel map $\psi \colon X' \to \erg$,
 and
 \item a Borel map $\xi \colon \erg \to \Prob(X')$,
 \end{itemize}
 such that
 \begin{enumerate}
 \item \label{ErgDecompThm-atom}
 $\xi(\omega) \bigl( \psi^{-1}(\omega) \bigr) = 1$ for each $\omega \in
\erg$,
 \item $\mu = \int_{\erg} \xi(\omega) \, d\nu(\omega)$,
 and
 \item For each $\omega \in \erg$, $\xi(\omega)$ is quasi-invariant and ergodic.
 \end{enumerate}
 \end{thm}

\begin{rem} \label{UnionOfProds}
 To simplify the notation in the conclusion of Thm.~\ref{ErgDecompThm}, we
will often assume that the space~$X$ can be written as a Cartesian product
$X = \erg \times S$, such that
 \begin{itemize}
 \item $\mu$ is the product measure on $X = \erg \times S$,
 and
 \item $\psi(\omega, s) = \omega$, for a.e.\ $(\omega,s) \in X$.
 \end{itemize}
 (For example, we make this assumption in the statement of
Prop.~\ref{BorelMaps}.)
 The following decomposition theorem of V.~A.~Rohlin \cite{Rohlin}
asserts that (up to isomorphism) the general case is a countable union of
examples of this type.
 \end{rem}

\begin{prop}[Rohlin] \label{RohlinDecomp}
 Assume the notation of Theorem~\ref{ErgDecompThm}. There is a partition
of~$\Omega$ into countably many Borel sets $\Omega_1,\Omega_2,\ldots$
{\rm(}some of these sets may be empty{\rm)}, such that, for each~$k$,
there exist
 \begin{enumerate}
 \item a conull subset $X'_k$ of $\psi^{-1}(\Omega_k)$,
 \item a standard Borel probability space $(S_k,\mu_k)$,
 \item a $(\psi_*\mu \times \mu_k)$-conull subset $(\erg_k \times
S_k)'$ of $\erg_k \times S_k$,
 and
 \item a measure-class-preserving Borel isomorphism $\theta \colon X'_k
\to (\erg_k \times S_k)'$,
 \end{enumerate}
 such that $\psi(x) = \pi_1 \bigl( \theta(x) \bigr)$
for each $x \in X'_k$, where $\pi_1(\omega,s) = \omega$.
 \end{prop}

\begin{proof}
 Say that two standard Borel probability spaces $(S_1,\mu_1)$ and
$(S_2,\mu_2)$ are \emph{of the same type} if there exists a
measure-class-preserving Borel isomorphism from a conull subset of~$S_1$
onto a conull subset of~$S_2$ \cite[pp.~10--11]{Rohlin}. This is obviously
an equivalence relation. It has only countably many equivalence classes
\cite[p.~18]{Rohlin}. Thus, there is a decomposition $\erg = \erg_1 \cup
\erg_2 \cup \cdots$, such that if $\omega$ and~$\omega'$ belong to the
same~$\Omega_k$, then the fibers $\bigl( \psi^{-1}(\omega), \xi(\omega)
\bigr)$ and $\bigl( \psi^{-1}(\omega'), \xi(\omega') \bigr)$ are of the
same type \cite[(I), p.~41]{Rohlin}. This implies that, modulo sets of
measure~$0$, $\psi^{-1}(\erg_k)$ is Borel isomorphic to the Cartesian
product $\erg_k \times \psi^{-1}(\omega)$, for any $\omega \in \erg_k$
\cite[p.~42]{Rohlin}.
 \end{proof}

\begin{prop} \label{BorelMaps}
 Let
 \begin{itemize}
 \item $G$ be a second countable, locally compact group,
 \item $(\erg,\nu)$ and $(S,\mu)$ be standard Borel probability
spaces,
 \item $\erg'$ be a conull Borel subset of~$\erg$,
 and
 \item $\rho \colon G \times (\erg \times S) \to \erg \times S$ be a
Borel action of~$G$ on $\erg \times S$, such that, for each $\omega \in
\erg'$, the probability measure $\mu_\omega$ on $\{\omega\} \times S$
{\rm(}induced by the natural isomorphism with~$S${\rm)} is
quasi-invariant.
 \end{itemize}
  Then:
 \begin{enumerate}
 \item  \label{BorelMaps-intoA}
 There are Borel maps $\rho_S \colon G \times
\erg' \times S \to S$ and $\rho_{\Aut} \colon G \times \erg' \to \A$,
defined by
 \begin{equation} \label{BorelMaps-RhoDefn}
   \mbox{$\rho(g,\omega,s) = \bigl( \omega, \rho_S(g,\omega,s) \bigr)
 = \bigl( \omega, \rho_{\Aut}(g,\omega)(s)  \bigr)$ for a.e.\ $s \in S$}
 \end{equation}
 for each $(g,\omega) \in G \times \erg'$.
 \item \label{BorelMaps-RN}
 There is a Borel function $D \colon G \times \erg \times S \to \real^{\ge0}$, such that
 $$\int_S D(g, \omega,s) f(s) \, d\mu(s) = \int_S f \bigl(
\rho_S(g,\omega,s) \bigr)\, d\mu(s) $$
 for every $g \in G$, $\omega \in \erg$, and $f \in \F(S,
\real^{\ge0})$.
 \end{enumerate}
 \end{prop}

\begin{rem}
The function $D$ in conclusion $2$ above is a \emph{fiberwise}
Radon-Nikodym derivative.  In particular, in the case when
$\Omega$ is a one point set, conclusion $2$  implies that the Radon-Nikodym
derivative is a Borel function on $G{\times}S$. \end{rem}

\begin{proof}
 Equation~\pref{BorelMaps-RhoDefn} determines a well-defined function
$\rho_{\Aut}  \colon G \times \erg' \to \F(S,S)$. Because $G$ is a group,
we know that $\rho_{\Aut}(g,\omega)$ is (essentially) a Borel
automorphism. Because $\mu_\omega$ is quasi-invariant, we see that
$\rho_{\Aut}(g,\omega)$ is measure-class preserving. Thus,
$\rho_{\Aut}(g,\omega) \in \A$, so $\rho_{\Aut}$ is actually a map
into $\A$. Thus, in order to complete the proof of
\pref{BorelMaps-intoA}, it only remains to show that $\rho_{\Aut}$ is
Borel. For this, we will use the conclusion of \pref{BorelMaps-RN}, so
let us establish the latter.

\pref{BorelMaps-RN}
 Let $\{A_n\}$ be a countable, dense subset of the Borel algebra $\B(S)$,
and define
 \begin{align*}
 D^\Delta &= \bigcap_{n \in \natural} \bigset{
 (g, \omega,f)
 }{
 \mu \bigl( \rho_{\Aut}(g, \omega)(A_n) \bigr)
 = \int_{A_n} f \, d\mu} \\
 &\subset G \times \erg' \times \F(S,\real^{\ge0})
 . \end{align*}
For each $g \in G$ and $\omega \in \erg'$, the fiber
 $\{\, s \in S \mid (g,\omega,s) \in D^\Delta \,\}$
 of~$D^\Delta$ consists precisely of the Radon-Nikodym derivative of the
transformation $\rho_{\Aut}(g, \omega)$. Therefore, $D^\Delta$ is the graph of
a function $\check D \colon G \times \erg' \to \F(S,\real^{\ge0})$, and
$\check D(g, \omega)$ is the Radon-Nikodym derivative of $\rho_{\Aut}(g,
\omega)$.

 Note:
 \begin{itemize}
 \item Because the Borel map $(g,\omega,s) \mapsto \bigl( g,\omega,
\rho_S(g,\omega,s) \bigr)$ is injective, we know that it maps $G \times
\erg' \times A_n$ (or any other Borel set) to a Borel subset of $G \times
\erg' \times S$. So Fubini's Theorem implies that $\mu \Bigl(
\rho_S \bigl( (g, \omega) \times A_n \bigr) \Bigr)$ is a Borel function of
$(g,\omega)$.
 \item Lemma~\ref{IntIsCont} implies that $\int_{A_n} f \, d\mu$ is a
Borel function of~$f$.
 \end{itemize}
 Therefore $D^\Delta$ is a Borel set, so the corresponding function
$\check D$ is a Borel function. Then the desired Borel function $D \colon
(G \times \erg') \times S \to \real^{\ge0}$ is obtained by applying
Lem.~\ref{DefdFiberwise}.

\pref{BorelMaps-intoA}
 Given $L^2$ functions $f,g \colon S \to \real$, and any $\epsilon > 0$,
define
 $$ \Theta \colon G \times \erg' \times S \to \real^{\ge0}$$
 by
 \begin{align*}
 \Theta(g,\omega,s)
 &= \bigl| T_{\rho_{\Aut}(g,\omega)} f(s) - g(s) \bigr|^2 \\
 &= \bigl| D(g,\omega,s)^{1/2} \, f\bigl( \rho_{\Aut}(g,\omega)^{-1}(s)
\bigr) - g(s) \bigr|^2 ,
 \end{align*}
 where $D(g,\omega,s)$ is given by \fullref{BorelMaps}{RN}.
 Then $\Theta$ is a Borel function on $G \times \erg' \times S$, so
 Fubini's Theorem implies that
 $$ \bigset{ (g, \omega) \in G \times \erg' }{ \int_S \Theta(g,\omega,s) <
\epsilon^2} $$
 is Borel. In other words, if $\mathcal{U}(f,g,\epsilon)$ is any basic
open set in $\Unitary$, then $\rho_{\Aut}^{-1} \bigl(
\mathcal{U}(f,g,\epsilon) \bigr)$ is Borel. So $\rho_{\Aut}$ is Borel.
 \end{proof}

\begin{cor} \label{BorelIntoF}
 Let
 \begin{itemize}
 \item $G$, $\erg$, $S$, $\erg'$, and $\rho$ be as in
Prop.~\ref{BorelMaps},
 and
 \item $X$ be a complete, separable metric space.
 \end{itemize}
 Then the action~$\rho$ induces a Borel map $\rho_{\F} \colon G \times
\erg' \times \F(S,X) \to \F(S,X)$, defined by
 \begin{equation}  \label{BorelIntoF-eq}
 \mbox{$\rho_{\F}(g,\omega,f)(s) = f \bigl( \rho_S(g,\omega,s) \bigr)$
 \quad for a.e.~$s$} .
 \end{equation}
 \end{cor}

\begin{proof}
 The existence of an abstract function~$\rho_{\F}$ satisfying
\pref{BorelIntoF-eq} is not an issue. Because $\rho_{\Aut}$ is Borel
\fullsee{BorelMaps}{intoA}, and the natural action of $\A$ on $\F(S,X)$ is
continuous \see{Autmu(S)action}, we know that $\rho_{\F}$ is Borel.
 \end{proof}

\subsection{Borel cocycles}

We assume the basic theory of Borel cocycles, as in
\cite[\S7.2]{MargulisBook} or \cite[\S4.2]{ZimmerBook}.

\begin{defn}[{\cite[Defns.~4.2.1 and 4.2.2]{ZimmerBook}}] \label{CocDefn}
 Suppose $H$~is a topological group, and $\rho \colon G \times X \to X$
is a Borel action of a locally compact group~$G$ on a standard Borel
probability space~$X$ with quasi-invariant measure~$\mu$.
 \begin{enumerate}
 \item A Borel function $\alpha \colon G \times X \to H$ is a Borel
\emph{cocycle} (for the action~$\rho$) if, for all $g_1,g_2 \in G$, we
have
 \begin{equation} \label{CocDefn-CocEqn}
 \mbox{$\alpha(g_1 g_2, x) = \alpha \bigl( g_1, \rho(g_2 , x) \bigr) \,
\alpha(g_2, x)$
 \quad for a.e.\ $x \in X$}.
 \end{equation}
 \item The cocycle $\alpha$ is \emph{strict} if the equality in
\pref{CocDefn-CocEqn} holds for \emph{all} $x \in X$, not merely almost
all.
 \item Two Borel cocycles $\alpha,\beta \colon G \times X \to H$ are
\emph{cohomologous} if there is a Borel function $\phi \colon X \to H$,
such that, for each $g \in G$, we have
 $$ \mbox{$\beta(g,x) = \phi \bigl( \rho(g,x) \bigr) \, \alpha(g,x) \,
\phi(x)^{-1}$
 \quad for a.e.\ $x \in X$} .$$
 This is an equivalence relation.
 \end{enumerate}
 \end{defn}

\begin{rem} \label{CanBeStrict}
 There is usually no harm in assuming that a Borel cocycle is strict,
because any Borel cocycle is
equal a.e.\ to a strict Borel cocycle \cite[Thm.~B.9, p.~200]{ZimmerBook}.
More precisely, if $\alpha$ is a Borel cocycle, then there is a strict
Borel cocycle $\alpha'$, such that, for each $g \in G$, $\alpha(g,x) =
\alpha'(g,x)$ for a.e.\ $x \in X$.
 \end{rem}

\begin{lem} \label{PointwiseOps}
 Let
 \begin{itemize}
 \item $(S,\mu)$ be a standard Borel probability space.
 \item $X$, $Y$, and~$Z$ be complete, separable, locally compact
metric spaces,
 \item $\tau \colon X \times Y \to Z$ be a continuous function.
 \end{itemize}
 Then the induced map $\tau^{\F} \colon \F(S,X) \times \F(S,Y) \to
\F(S,Z)$, defined by
 $$ \tau^{\F}(\phi,\psi)(s) = \tau \bigl( \phi(s), \psi(s) \bigr), $$
 is continuous.
 \end{lem}

\begin{proof}
 Given sequences $\phi_n \to \phi \in \F(S,X)$, $\psi_n \to \psi \in
\F(S,Y)$, and $\epsilon > 0$, Lusin's Theorem gives us a compact
subset~$K$ of~$S$, such that $\phi$ and~$\psi$ are continuous on~$K$, and
$\mu(K) > 1 - \epsilon$. Because $X$ and~$Y$ are locally compact, we may
let $K_X$ and $K_Y$ be compact neighborhoods of $\phi(K)$ and $\psi(K)$
in $X$ and~$Y$, respectively.
 Because $\tau$ is uniformly continuous on $K_X \times K_Y$, there is
some $\delta > 0$, such that
 \begin{align*}
 &\limsup_{n \to \infty}  \mu \bigset{ s \in S }{
 d \Bigl( \tau \bigl( \phi_n(s), \psi_n(s) \bigr), \tau \bigl( \phi(s),
\psi(s) \bigr) \Bigr) > \epsilon } \\
 &\kern3em\le
 \limsup_{n \to \infty} \mu \{\, s \in S \mid \mbox{$\phi_n(s) \notin
K_X$ or $\psi_n(s) \notin K_Y$} \,\} \\
 &\kern6em + \limsup_{n \to \infty} \mu \{\, s \in S \mid d
\bigl( \phi_n(s), \phi(s) \bigr) > \delta \,\} \\
 &\kern6em + \limsup_{n \to \infty} \mu \{\, s \in S \mid d
\bigl( \psi_n(s), \psi(s) \bigr) > \delta \,\} \\
 &\kern3em \le 2 \epsilon + 0 + 0
 . \end{align*}
 Because $\epsilon > 0$ is arbitrary, we conclude that
 $\tau^{\F}(\phi_n,\psi_n) \to \tau^{\F}(\phi,\psi)$.
 \end{proof}

\begin{cor}[{cf.\ \cite[Rmk.\ after Lem.~7.2.1, p.~217]{MargulisBook}}]
\label{TwistisBorel}
 Let
 \begin{itemize}
 \item $G$, $\erg$, $S$, $\erg'$, and $\rho$ be as in
Prop.~\ref{BorelMaps},
 \item $H$ be a locally compact, second countable group,
 \item $X$ be a complete, separable metric space,
 \item $\tau \colon H \times X \to X$ be a continuous action of~$H$
on~$X$,
 and
 \item $\alpha \colon G \times (\erg \times S) \to H$ be a Borel cocycle.
 \end{itemize}
 Then the function $\rho^{\alpha,\tau} \colon G \times \erg' \times
\F(S,X) \to \F(S,X)$, defined by
 $$ \rho^{\alpha,\tau}(g,\omega,\phi)(s) =
 \tau \Bigl( \alpha(g,\omega, s)^{-1}, \phi \bigl( \rho_S(g,\omega,s)
\bigr) \Bigr) ,$$
 is Borel.
 \end{cor}

\begin{proof}
 Because the map $\rho_{\F}$ of Cor.~\ref{BorelIntoF} is Borel, the
action~$\rho$ does not affect the measurability of~$\rho^{\alpha,\tau}$,
so it may be ignored. Furthermore, by replacing $\erg$ with $G \times
\erg'$, we may assume that $G$~is trivial and $\erg' = \erg$; in
particular, $G$ may be ignored.  Thus,
 \begin{itemize}
 \item $\alpha$ is a Borel map from $\erg \times S$ to $H$,
 and
 \item $\rho^{\alpha,\tau} \colon \erg  \times \F(S,X) \to \F(S,X)$ is
defined by
 $$ \rho^{\alpha,\tau}(\omega,\phi)(s) =
 \tau \bigl( \alpha(\omega,s)^{-1}, \phi(s) \bigr) .$$
 \end{itemize}
 Then, because $\check\alpha \colon \erg \to \F(S,X)$ and $\tau^{\F}$
are Borel (see \fullref{FiberFunc}{Borel} and~\ref{PointwiseOps}), we
conclude that $\rho^{\alpha,\tau}$ is a composition of Borel functions.
Therefore, it is Borel.
 \end{proof}

Although we do not need the following result in this paper, we include
the proof to provide a convenient reference.

\begin{cor}[{\cite[Rmk.\ after Lem.~7.2.1, p.~217]{MargulisBook}}]
\label{TwistContinuous}
 Let
 \begin{itemize}
 \item $G$ and $H$ be second countable, locally compact groups,
 \item $(S,\mu)$ be a standard Borel probability space,
 \item $\rho \colon G \times S \to S$ be a Borel action of~$G$ on~$S$,
such that  $\mu$ is quasi-invariant,
 \item $\alpha \colon G \times S \to H$ be a strict Borel cocycle,
 and
 \item $\tau \colon H \times X \to X$ be a continuous action of~$H$
on~$X$.
 \end{itemize}
 Then
 \begin{enumerate}
 \item the function $\check\alpha \colon G \to \F(S,H)$, induced
by~$\alpha$ \see{FiberFunc}, is continuous,
 and
 \item the $\alpha$-twisted action $\zeta_{\rho,\tau,\alpha}$ of~$G$ on
$\F(S,H)$ is continuous {\rm(}see Defn.~\fullref{TwistDefn}{action}{\rm)}.
 \end{enumerate}
 \end{cor}

\begin{proof}
  Let $\triv \colon G \times S \to H$ be the trivial cocycle, defined
by $\triv(g,s) = e$ (the identity element of~$H$). For convenience, let
$\F = \F(S,H)$.

\setcounter{steppf}{0}

\begin{steppf} \label{TwistContinuous-translate}
 The action $\rho^{\triv,\alpha}$ is continuous.
 \end{steppf}
 From Prop.~\fullref{BorelMaps}{intoA} (with $\erg$ consisting of a single
point), we know that $\rho$ induces a Borel function $\rho_{\Aut} \colon
G \to \A$. This function is a homomorphism, and any measurable
homomorphism into a second countable topological group is continuous
\cite[Thm.~B.3, p.~198]{ZimmerBook}, so we conclude that $\rho_{\Aut}$ is
continuous. Because the action of $\A$ on $\F$ is continuous
\see{Autmu(S)action}, we conclude that $\rho^{\triv,\alpha}$ is
continuous.

\begin{steppf} \label{TwistContinuous-cocycle}
 $\check\alpha$ is continuous.
 \end{steppf}
 Note that $\F$ is a (second countable) topological group under
pointwise multiplication \cf{PointwiseOps}. Furthermore,
$\rho^{\triv,\alpha}$ is a continuous action of~$G$ on~$\F$ by
automorphisms, so we may form the semidirect product $G \ltimes \F$. The
cocycle identity implies that the Borel function $\alpha^{G \ltimes \F}
\colon G \to G \ltimes \F$, defined by $\alpha^{G \ltimes \F}(g,f)
= \bigl( g, \check\alpha(g) \bigr)$ is a homomorphism. Because (as noted
above) measurable homomorphisms are continuous, we conclude that
$\check\alpha$ is continuous.

\begin{steppf}
 $\rho^{\tau,\alpha}$ is continuous.
 \end{steppf}
 Because $\rho^{\tau,\alpha}(g,f) = \tau^{\F} \bigl( \check\alpha(g),
\rho^{\triv,\alpha}(g,f) \bigr)$, this conclusion can be obtained by
combining Steps~\ref{TwistContinuous-translate}
and~\ref{TwistContinuous-cocycle} with Lem.~\ref{PointwiseOps}.
 \end{proof}

\section{Restricting cocycles to ergodic components}
\label{RestrictCocycleSect}

This section uses the von Neumann Selection Theorem \pref{vonNeumann} to
obtain information about a cocycle from its restrictions to ergodic
components. The main result is \pref{AnalErgComp}; the others are
corollaries.

\begin{notn} \label{ResCocNotn}
 Let
 \begin{itemize}
 \item $G$, $\erg$, $S$, $\erg'$ and~$\rho$ be as in
Prop.~\ref{BorelMaps} (or, equivalently, as in Thm.~\ref{AnalErgComp}
below),
 \item $H$ be a locally compact, second countable group,
 and
 \item $\alpha \colon G \times (\erg \times S) \to H$ be a strict Borel
cocycle.
 \end{itemize}
 For each $\omega \in \erg'$, define
 $\rho_\omega \colon G \times S \to S$
 and
 $\alpha_\omega \colon G \times S \to H$
 by
 $$ \mbox{$\rho_\omega(g,s) = \rho_S(g,\omega,s)$
 \quad and \quad
 $\alpha_\omega(g,s) = \alpha(g,\omega,s)$}
 ,$$
 where $\pi_S(*,s) = s$.
 \end{notn}

\begin{defn}[{\cite[Lem.~3.1]{Ramsay}, \cite[Defn.~3.1]{ZimmerExtn}}]
 Suppose $G$ is a locally compact second countable group, and $(S,\mu)$ is
a standard Borel probability space.
 A Borel map $\rho \colon G \times S \to S$ is a \emph{near action}
of~$G$ on~$S$ if
 \begin{itemize}
 \item for all $g_1,g_2 \in G$, we have
 $\rho(g_1 g_2 , s) = \rho \bigl( g_1, \rho(g_2, s) \bigr)$
 for a.e.\ $s \in S$,
 \item $\rho(e,s) = s$ for a.e.\ $s \in S$,
 and
 \item each $g \in G$ preserves the measure class of~$\mu$.
 \end{itemize}
 Note that the definition of a cocycle \pref{CocDefn} can be applied to
near actions, not only actions.
 \end{defn}

Let us record the following elementary observation.
Part~\pref{ResCocLem-Invt} follows easily from the assumption that
$\mu_\omega$ is quasi-invariant. The other two parts are consequences of
the first.

\begin{lem} \label{ResCocLem}
 Assume the setting of Notation~\ref{ResCocNotn}, and let $\omega \in
\erg'$. Then:
 \begin{enumerate}
 \item \label{ResCocLem-Invt}
 $\rho(\omega,s) \in \{\omega\} \times S$, for a.e.\ $s \in S$,
 \item $\rho_\omega$ is a near action of~$G$ on~$S$,
 and
 \item \label{ResCocLem-coc}
 $\alpha_\omega$ is a Borel cocycle for~$\rho_\omega$.
 \end{enumerate}
 \end{lem}

\begin{thm} \label{AnalErgComp}
 Let
 \begin{itemize}
 \item $G$ and~$H$ be second countable, locally compact groups,
 \item $(\erg,\nu)$ and $(S,\mu)$ be standard Borel probability
spaces,
 \item $\erg'$ be a conull subset of~$\erg$,
 \item $\rho \colon G \times (\erg \times S) \to \erg \times S$ be a
Borel action of~$G$ on $\erg \times S$, such that, for each $\omega' \in
\erg$, the probability measure $\mu_\omega$ on $\{\omega\} \times S$
{\rm(}induced by the natural isomorphism with~$S${\rm)} is
quasi-invariant,
 \item $\alpha \colon G \times (\erg \times S) \to H$ be a strict Borel
cocycle,
 \item $\func$ be an analytic subset of $\erg \times \F(G \times
S,H)$,
 and
 \item $\erg_\func = \{\, \omega \in \erg' \mid \mbox{$\alpha_\omega$
is cohomologous to a cocycle in $\func_\omega$} \,\}$,
 where
 $$\func_\omega = \{\, f \in \F(G \times
S,H) \mid (\omega, f) \in \func\,\} .$$
 \end{itemize}
  Then:
 \begin{enumerate}
 \item \label{AnalErgComp-anal}
 $\erg_\func$ is analytic,
 and
 \item \label{AnalErgComp-coho}
 $\alpha$ is cohomologous to a Borel cocyle
$\beta \colon G \times (\erg \times S) \to H$, such that $\beta_\omega \in
\func_\omega$, for a.e.\ $\omega \in \erg_{\func}$.
 \end{enumerate}
 \end{thm}

\begin{proof}
 The cocycle~$\alpha$ defines a Borel map $\check\alpha \colon \erg \to
\F(G \times S,H)$ \fullsee{FiberFunc}{Borel}.  Define
 $$ \cobound \colon \erg' \times \F(G \times S,H) \times \F(S,H) \to \F(G
\times S,H)$$
 by
 \begin{equation} \label{AnalErgComp-delta}
 \cobound(\omega, \phi, f)(g,s) = f \bigl( \rho_\omega(g,s) \bigr)
\, \phi(g,s) \, f(s)^{-1} .
 \end{equation}

We claim that $\cobound$ is Borel.
 \begin{itemize}
 \item We know that the map $\rho_{\F}$ (defined in Cor.~\ref{BorelIntoF})
is Borel. This induces a Borel map $\check\rho_{\F} \colon \erg' \times
\F(S,H) \to \F \bigl( G, \F(S,H) \bigr)$ \see{FiberFunc}. From
Cor.~\ref{F(F)}, we see that we may think of this as a map into $\F(G
\times S, H)$. Thus, the first factor on the right-hand side of
\pref{AnalErgComp-delta} represents a Borel function from $\erg' \times
\F(S,H)$ into  $\F(G \times S, H)$.
 \item The second factor on the right-hand side of
\pref{AnalErgComp-delta} represents the identity function on $\F(G \times
S, H)$, and the term $f(s)$ represents the inclusion of $\F(S,H)$ into
$\F(G \times S, H)$. These are obviously Borel maps into $\F(G \times S,
H)$.
 \end{itemize}
 Because pointwise multiplication and pointwise inversion are continuous
operations on $\F(G \times S, H)$ \see{PointwiseOps}, we conclude that
$\cobound$ is Borel, as claimed.

 Therefore, the function
 $$\sigma \colon \erg' \times \F(S,H) \to \erg' \times \F(G \times S,H)
,$$
 defined by
 $$ \sigma(\omega, f) = \Bigl( \omega, \cobound \bigl( \omega,
\check\alpha(\omega), f \bigr) \Bigr) ,$$
 is Borel, so $\sigma^{-1}(\func)$ is analytic. Then, because $\erg_\func$
is the projection of~$\sigma^{-1}(\func)$ to~$\erg'$, we conclude that
$\erg_\func$ is analytic. This establishes \pref{AnalErgComp-anal}.

The von Neumann Selection Theorem \pref{vonNeumann} implies that there is
a Borel function $\hat\Phi \colon \erg' \to \F(S,H)$, such that
$\sigma\bigl( \omega, \hat\Phi(\omega) \bigr) \in \func$, for a.e.\
$\omega \in \erg_\func$.  Corresponding to~$\hat\Phi$, there is a Borel
function $\Phi \colon \erg' \times S \to H$ \see{DefdFiberwise}. Letting
 $$ \beta(g,\omega,s)
 = \cobound \bigl( \omega, \check\alpha(\omega), \hat\Phi(\omega) \bigr)
(g,s)
 = \Phi \bigl( \rho(g,\omega,s) \bigr) \, \alpha(g,\omega,s) \,
\Phi(\omega, s)^{-1} ,$$
 we obtain \pref{AnalErgComp-coho}.
 \end{proof}

\begin{cor} \label{CohoErg}
 Let
 \begin{itemize}
 \item $G$ and $H$ be locally compact, second countable groups,
 \item $(X,\mu)$ be a standard Borel probability space,
 \item $\rho \colon G \times X \to X$ be a Borel action, such that $\mu$
is quasi-invariant,
 \item $\psi \colon X' \to \erg$ be the corresponding ergodic
decomposition \see{ErgDecompThm},
 and
 \item $\alpha,\beta \colon G \times X \to H$ be Borel cocycles.
 \end{itemize}
 For each $\omega \in \erg$, let $\alpha_\omega$ and $\beta_\omega$ be
the restrictions of~$\alpha$ and~$\beta$ to $G \times \psi^{-1}(\omega)$.
 Then:
 \begin{enumerate}
 \item \label{CohoErg-coc}
 There is a conull Borel subset $\erg'$ of~$\erg$, such that, for each
$\omega \in \erg'$, the maps $\alpha_\omega$ and $\beta_\omega$ are Borel
cocycles.
 \item \label{CohoErg-anal}
 $\{\, \omega \in \erg' \mid \mbox{$\alpha_\omega$ is cohomologous
to~$\beta_\omega$}\,\}$ is an analytic subset of~$\erg$.
 \item \label{CohoErg-trivial}
 If $\alpha_\omega$ is cohomologous to~$\beta_\omega$, for a.e.\ $\omega
\in \erg'$, then $\alpha$ is cohomologous to~$\beta$.
 \end{enumerate}
 \end{cor}

\begin{proof}
 From Prop.~\ref{RohlinDecomp} (and Rem.~\ref{UnionOfProds}), we may
assume the notation of Thm.~\ref{AnalErgComp}. By changing $\alpha$
and~$\beta$ on a null set, we may assume these cocycles are strict
\see{CanBeStrict}.

Conclusion \pref{CohoErg-coc} is immediate from
Lem.~\fullref{ResCocLem}{coc}.

 Recall that $\beta$ induces a Borel function $\check\beta \colon \erg'
\to \F(G \times S, H)$, defined by $\check\beta(\omega) = \beta_\omega$
\fullsee{FiberFunc}{Borel}.
 Let
 $$\func = \{\, (\omega, \check\beta(\omega) \mid \omega \in \erg' \,\}
\subset \erg' \times \F(G \times S, H) .$$
 Because $\func$ is an analytic subset (in fact, it is closed),
\pref{CohoErg-anal} is immediate from \fullref{AnalErgComp}{anal}.

Assume, now, that $\alpha_\omega$ is cohomologous to~$\beta_\omega$, for a.e.\ $\omega
\in \erg'$. From \fullref{AnalErgComp}{coho} and Fubini's
Theorem, we conclude that $\alpha$ is cohomologous to a
cocycle~$\tilde\alpha$, such that for a.e.\ $g \in G$,
 \begin{equation} \label{CohoErg-AlmostBeta}
 \mbox{for a.e.~$x \in \erg' \times S$, we have $\tilde\alpha(g,x) =
\beta(g,x)$.}
 \end{equation}
 From the cocycle identity, one easily concludes that
\pref{CohoErg-AlmostBeta} must hold for every $g \in G$, not merely
almost every~$g$. Therefore $\tilde\alpha$ is (obviously) cohomologous
to~$\beta$. By transitivity, then $\alpha$ is also cohomologous
to~$\beta$; this establishes~\pref{CohoErg-trivial}.
 \end{proof}

\begin{defn}
 Suppose $\rho \colon G \times X \to X$ is a Borel action with
quasi-invariant measure, and $H$~is a locally compact second countable
group.
 \begin{enumerate}
 \item The \emph{trivial cocycle} $\triv_{G \times X} \colon G \times X
\to H$ is defined by $\triv_{G\times X}(g,x) = e$.
 \item A Borel cocycle $\alpha \colon G \times X \to H$ is a
\emph{coboundary} if it is cohomologous to the trivial cocycle.
 \end{enumerate}
 \end{defn}

\begin{cor} \label{TrivialCocycle}
 Let $G$, $H$, $X$, $\rho$, $\psi$, $\Omega$, $\alpha$, $\alpha_\omega$,
and~$\erg'$ be as in Cor.~\ref{CohoErg}. Then:
 \begin{enumerate}
 \item \label{TrivialCocycle-anal}
 $\{\, \omega \in \erg' \mid \mbox{$\alpha_\omega$ is a coboundary}\,\}$
is an analytic subset of~$\erg$.
 \item \label{TrivialCocycle-trivial}
 If $\alpha_\omega$ is a coboundary, for a.e.\ $\omega \in \erg'$, then
$\alpha$ is a coboundary.
 \end{enumerate}
 \end{cor}

\begin{proof}
 Let $\beta$ be the trivial cocycle, and apply Cor.~\ref{CohoErg}.
 \end{proof}

\begin{defn} \label{HomoCocDefn}
 Recall that a Borel cocycle $\alpha \colon G \times S \to H$ is a
\emph{constant} (or \emph{homomorphism}) cocycle if $\alpha(g,s)$ is
essentially independent of~$s$, for each $g \in G$.
 \end{defn}

\begin{cor} \label{CohoHomo}
 Let $G$, $H$, $X$, $\rho$, $\psi$, $\Omega$, $\alpha$, $\alpha_\omega$,
and~$\erg'$ be as in Cor.~\ref{CohoErg}.
 If $\alpha_\omega$ is cohomologous to a constant cocycle, for a.e.\
$\omega \in \erg'$, then $\alpha$ is cohomologous to a constant cocycle.
 \end{cor}

\begin{proof}
 Similar to the proof of Cor.~\fullref{CohoErg}{trivial}, but with
 $\func = \erg' \times \Const$, where
 $$ \Const = \{\, f \colon G \times S \to H \mid \mbox{$f(g,s)$ is
essentially independent of~$s$} \,\} .$$
 \end{proof}

\begin{notn}
 For any locally compact, second countable group~$H$,
we use $\cpct(H)$ to denote the set of compact subgroups of~$H$. It is
well known that this is a complete, separable metric space, under the
Hausdorff metric
 $$ d(K_1,K_2) = \max_{k_1 \in K_1} \dist(k_1,K_2) + \max_{k_2 \in K_2}
\dist(K_1,k_2) ,$$
 where $\dist$ is any metric on~$H$.
 \end{notn}

\begin{cor} \label{CohoCpctImg}
 Let $G$, $H$, $X$, $\rho$, $\psi$, $\Omega$, $\alpha$, $\alpha_\omega$,
and~$\erg'$ be as in Cor.~\ref{CohoErg}.

If, for a.e.\ $\omega \in \erg'$, the cocycle $\alpha_\omega$ is
cohomologous to a cocycle whose essential range is contained in a compact
subgroup of~$H$, then there are
 \begin{itemize}
 \item a Borel function $\kappa \colon \erg' \to \cpct(H)$,
 and
 \item a Borel cocycle $\beta$ that is cohomologous to~$\alpha$,
 \end{itemize}
 such that the essential range of~$\beta_\omega$ is
contained in~$\kappa(\omega)$, for a.e.\ $\omega \in \erg'$.
 \end{cor}

\begin{proof}
 As in the proof of Cor.~\ref{CohoErg}, we assume the notation of
Thm.~\ref{AnalErgComp}, and we assume the cocycle~$\alpha$ is strict.
 Let
 $$ \func^+ = \{\, (f,K) \in \F(G \times S, H) \times \cpct(H) \mid
\EssRg(f) \subset K \,\} ,$$
 and let $\func$ be the projection of~$\func^+$ to $\F(G \times S, H)$.
Then $\func^+$ is closed, so $\func$ is analytic. Applying
Thm.~\fullref{AnalErgComp}{coho} yields a Borel cocycle~$\beta$,
cohomologous to~$\alpha$, such that $\beta_\omega \in F$, for a.e.\
$\omega \in \erg'$. Now the Borel function $\kappa$ is obtained from
the von Neumann Selection Theorem \pref{vonNeumann}.
 \end{proof}

\begin{cor} \label{HomoModCpct}
 Let $G$, $H$, $X$, $\rho$, $\psi$,
$\Omega$, $\alpha$, $\alpha_\omega$, and~$\erg'$ be as in
Cor.~\ref{CohoErg}.

 If, for a.e.\ $\omega \in \erg'$, there are a compact subgroup~$K_\omega$
of~$H$ and a Borel cocycle~$\beta_\omega$, cohomologous
to~$\alpha_\omega$, such that
 \begin{enumerate} \renewcommand{\theenumi}{\alph{enumi}}
 \item \label{HomoModCpct-range}
 the essential range of~$\beta_\omega$ is contained in the normalizer
$N_H(K_\omega)$,
 and
 \item \label{HomoModCpct-homo}
 the induced cocycle $\overline{\beta_\omega} \colon G \times S \to
N_H(K_\omega)/K_\omega$ is a homomorphism cocycle,
 \end{enumerate}
 then there are
 \begin{enumerate}
 \item a Borel cocycle~$\beta$, cohomologous to $\alpha$,
 and
 \item a Borel function $\kappa \colon \erg \to \cpct(H)$,
 \end{enumerate}
 such that \pref{HomoModCpct-range} and~\pref{HomoModCpct-homo} hold with
$\beta_\omega = \check\beta(\omega)$ and $K_\omega = \kappa(\omega)$, for
a.e.\ $\omega \in \erg'$.
 \end{cor}

\begin{proof}
  Let
 $$ \func^+ = \bigset{ (f,K) \in \F \bigl( G, \F(S, H) \bigr) \times
\cpct(H) }{\
 \begin{matrix}
 \mbox{for a.e.\ $g \in G$, $ \exists h \in N_H(K)$, } \\
 \mbox{such that $\EssRg\bigl( f(g) \bigr) \subset h \, K$}
 \end{matrix}
 } .$$
 Then $\func^+$ is closed, so the proof is similar to that of
Cor.~\ref{CohoCpctImg}. (Recall that $\F \bigl( G, \F(S, H) \bigr)$ is
naturally homeomorphic to $\F(G \times S, H)$ \see{F(F)}.)
 \end{proof}

\section{Superrigidity for non-ergodic actions}
\label{section:superrigid}

One main application of the results in this paper is to prove
general versions of superrigidity for cocycles for non-ergodic
actions of certain groups.  We now define this class of groups.  Let $I$ be a
finite index set and for each $i{\in}I$, we let $k_i$ be a local
field of characteristic zero and ${\mathbb G}_i$ be a connected
simply connected semisimple algebraic $k_i$-group. We first define
groups $G_i$, and then let $G=\prod_{i{\in}I}G_i$.  If $k_i$ is
non-Archimedean,
 $G_i={\mathbb G}_i(k_i)$ the
$k_i$-points of ${\mathbb G}_i$. If $k_i$ is Archimedean, then
$G_i$ is either ${\mathbb G}_i(k_i)$ or its topological universal
cover. (This makes sense, since when $\Ga_i$ is simply connected
and $k_i$ is Archimedean, ${\mathbb G}_i(k_i)$
 is topologically connected.) We assume that the $k_i$-rank of any
 simple factor of any ${\mathbb G}_i$ is at least two.

We will need one assumption on the cocycles we consider.

\begin{defn}
\label{defn:quasi-integrable} Let $D$ be a locally compact group,
$(S,\mu)$ a standard probability measure space on which $D$ acts
preserving $\mu$  and $H$ be a normed topological group. We call a
cocycle $\alpha\colon D{\times}S{\rightarrow}H$ over the $D$
action \emph{$D$-integrable} if for any compact subset
$M\subset{D}$, the function
$Q_{M,\alpha}(x)=\sup_{m{\in}M}\ln^+\|\alpha(m,x)\|$ is in
$L^1(S)$ (recall that $\ln^+ x = \max{(\ln x, 0)}$).
\end{defn}

Any continuous cocycle over a continuous action on a compact
topological space is automatically $D$-integrable.  We remark that
a cocycle over a cyclic group action is $D$-integrable if and only
if $\ln^+\|(\alpha(\pm{1},x)\|$ is in $L^1(S)$.

We first recall the superrigidity theorems from \cite{FisherMargulis} for ergodic actions.

\begin{thm}
\label{theorem:Gsuperrigidity} Let $G$ be as above, let $(S,\mu)$
be a standard probability measure space and let $H$ be the $k$
points of a $k$-algebraic group where $k$ is a local field of
characteristic $0$.  Assume $G$ acts ergodically on $S$ preserving
$\mu$. Let $\alpha\colon G{\times}S{\rightarrow}H$ be a
$G$-integrable Borel cocycle.  Then $\alpha$ is cohomologous to a
cocycle $\beta$ where $\beta(g,x)=\pi(g){c(g,x)}$. Here
${\pi\colon G{\rightarrow}H}$ is a continuous homomorphism and
$c\colon G{\times}S{\rightarrow}C$ is a cocycle taking values in a
compact group centralizing $\pi(G)$.
\end{thm}

\begin{thm}
\label{theorem:Gammasuperrigidity} Let $G,S,H$ and $\mu$ be as
Theorem \ref{theorem:Gsuperrigidity} and let $\G<G$ be a lattice.
Assume $\Gamma$ acts ergodically on $S$ preserving $\mu$. Assume
$\alpha\colon \Gamma{\times}S{\rightarrow}H$ is a
$\Gamma$-integrable, Borel cocycle. Then $\alpha$ is cohomologous
to a cocycle $\beta$ where
$\beta(\gamma,x)=\pi(\gamma){c(\gamma,x)}$. Here ${\pi\colon
G{\rightarrow}H}$ is a continuous homomorphism of $G$ and $c\colon
\Gamma{\times}X{\rightarrow}C$ is a cocycle taking values in a
compact group centralizing $\pi(G)$.
\end{thm}

To state non-ergodic versions of the above theorems, we will need
a Borel structure on the space of homomorphisms for $G$ to $H$.
Given $G$ as above and $H$ as in Theorem
\ref{theorem:Gsuperrigidity}, it is well-known that there are only
finitely many conjugacy classes of homomorphisms $\pi\colon
G{\rightarrow}H$.  We choose a set $\Pi=\{\pi_i\}$ of
representatives and endow it with the discrete topology, so as to
be able to consider measurable maps to $\Pi$.

Given a group $D$ acting on a standard probability measure space
$(X,\mu)$, we denote by $\erg$ the space of ergodic components of
the action and let $p\colon S{\rightarrow}\erg$ be the natural
projection.  We now state the general versions of the
superrigidity theorems above.

\begin{thm}
\label{theorem:Gsuperrigidityne} Let $G$ be as above, let
$(X,\mu)$ be a standard probability measure space and let $H$ be
the $k$ points of a $k$-algebraic group where $k$ is a local field
of characteristic $0$.  Assume $G$ acts on $S$ preserving $\mu$.
Let $\alpha\colon G{\times}S{\rightarrow}H$ be a $G$-integrable
Borel cocycle.  Then there exist measurable maps
$\pi\colon\erg{\rightarrow}\Pi,
\kappa\colon\erg{\rightarrow}\cpct(H)$ and $\phi\colon
S{\rightarrow}H$ with $\kappa(p(x)){\subset}Z_H(\pi(p(x))$ almost
everywhere such that
 $$\alpha(g,x)=\phi(gx)^{-1}\beta(g,x)\phi(x)$$
 where $\beta(g,x)=\pi(p(x))(g){c(g,x)}$. Here $c\colon
G{\times}S{\rightarrow}H$ is a measurable cocycle with
$c(g,x){\in}\kappa(p(x))$ almost everywhere.
\end{thm}

\begin{thm}
\label{theorem:Gammasuperrigidityne} Let $G,S,H$ and $\mu$ be as
Theorem \ref{theorem:Gsuperrigidityne} and let $\G<G$ be a
lattice. Assume $\Gamma$ acts on $X$ preserving $\mu$. Assume
$\alpha\colon\Gamma{\times}S{\rightarrow}H$ is a
$\Gamma$-integrable, Borel cocycle. Then there exist measurable
maps $\pi\colon\erg{\rightarrow}\Pi,
\kappa\colon\erg{\rightarrow}\cpct(H)$ and $\phi\colon
S{\rightarrow}H$ with $\kappa(p(x)){\subset}Z_H(\pi(p(x)))$ such
that $\alpha(\g,x)=\phi({\g}x)^{-1}\beta(\g,x)\phi(x)$ where
$\beta(\g,x)=\pi(p(x))(\g){c(\g,x)}$. Here
$c\colon\G{\times}S{\rightarrow}H$ is a measurable cocycle with
$c(\g,x){\in}\kappa(p(x))$ almost everywhere. \end{thm}

\begin{proof}[Proof of Theorems \ref{theorem:Gsuperrigidityne} and
\ref{theorem:Gammasuperrigidityne}] These are an immediate
consequence of Theorems \ref{theorem:Gsuperrigidity} and
\ref{theorem:Gammasuperrigidity}, Corollary \ref{HomoModCpct}, and
Proposition \ref{RohlinDecomp}.   Moreover, one can also prove
these results by using the proof of Theorems
\ref{theorem:Gsuperrigidity} and \ref{theorem:Gammasuperrigidity}
from \cite{FisherMargulis}, Proposition  \ref{RohlinDecomp} and
Corollary \ref{AlgHull} below.
\end{proof}

There are also versions of Theorems \ref{theorem:Gsuperrigidityne}
and \ref{theorem:Gammasuperrigidityne} which do not require that
we assume the cocycle is $G$-integrable and versions, in that
context, where the class of $G$ considered can be somewhat
broader, i.e. $G$ of rank at least $2$, with some/all simple
factors of rank $1$. To remove the $G$-integrability assumption
requires assumptions on the algebraic hull of the cocycle, while
weakening the rank assumption requires both assumptions on the
algebraic hull  and the assumption that the $G$ action on each
ergodic component of $(X,\mu)$ is weakly irreducible. These
assumptions are less natural in the non-ergodic setting, so we
leave it to the interested reader to formulate and prove such
results, using Theorems $3.6$ and $3.7$ of \cite{FisherMargulis}
in place of Theorems \ref{theorem:Gsuperrigidity} and
\ref{theorem:Gammasuperrigidity} above.

\section{Equivariant maps on ergodic components}

This section uses von Neumann Selection Theorem \pref{vonNeumann} to
prove that if almost every ergodic component of a $G$-action has a
fixed standard Borel $G$-space~$X$ as a measurable quotient, then
$X$~is a measurable quotient of the entire action. Actually, the
conclusion is proved in a more general setting that includes twisting by
cocycles \see{QuotErgComps}. This yields a corollary \pref{AlgHull} that
obtains information about a cocycle from the algebraic hulls of its
restrictions to ergodic components.

\begin{defn}[{\cite[\S3.2.0, pp.~216--217]{MargulisBook}}]
\label{TwistDefn}
 Suppose
 \begin{itemize}
 \item $G$ and $H$ are locally compact, second countable groups,
 \item $S$ and $X$ are Borel spaces,
 \item $\rho \colon G \times S \to S$ and $\tau \colon H \times X \to X$
are Borel actions,
 \item $\mu$ is a probability measure on~$S$,
 and
 \item $\alpha \colon G \times S \to H$ is a Borel cocycle.
 \end{itemize}
 Then:
 \begin{enumerate}
 \item \label{TwistDefn-action}
 We define an action $\zeta_{\rho,\tau,\alpha} \colon G \times \F(S,X)
\to \F(S,X)$ by
 $$ \zeta_{\rho,\tau,\alpha}(g,\phi)(s) = \tau \Bigl( \alpha(g,s), \phi
\bigl( \rho(g^{-1},s) \bigr) \Bigr) .$$
 We may refer to $\zeta_{\rho,\tau,\alpha}$ as the
\emph{$\alpha$-twisted action} of~$G$ on $\F(S,X)$. It is Borel
\see{TwistisBorel}.
 \item \label{TwistDefn-equi}
 A function $\phi \colon S \to X$ is \emph{essentially $(\rho,\tau,
\alpha)$-equivariant} if, for each $g \in G$, we have
 $$\mbox{$\phi \bigl( \rho(g,s) \bigr) = \tau \bigl( \alpha(g,s),
\phi(s) \bigr)$ for a.e.~$s \in S$}. $$
 In other words, a Borel function $\phi \colon S \to X$ is essentially
$(\rho,\tau, \alpha)$-equivariant if and only if it represents a fixed
point of the $\alpha$-twisted action of~$G$ on $\F(S,X)$.
 \end{enumerate}
 \end{defn}

\begin{prop} \label{AlphaEquiComps}
 Let
 \begin{itemize}
 \item $G$, $H$, $X$, $\rho$, $\psi$, $\alpha$, and~$\alpha_\omega$ be as
in Cor.~\ref{CohoErg},
 \item $\rho_\omega$ be the restriction of~$\rho$ to $G
\times \psi^{-1}(\omega)$, for each $\omega \in \erg$,
 \item $(Y,d)$ be a complete, separable metric space,
 \item $\tau \colon H \times Y \to Y$ be a continuous action of~$H$
on~$Y$,
 \item $\erg'$ be a conull Borel subset of~$\erg$, such that
$\alpha_\omega$ is a Borel cocycle, for each $\omega \in \erg'$,
 and
 \item $\erg^X = \bigset{ \omega \in \erg' }{
 \begin{matrix}
 \mbox{there is a Borel map $\phi_\omega \colon \psi^{-1}(\omega) \to Y$}
\\
 \mbox{that is essentially $(\rho_\omega, \tau,
\alpha_\omega)$-equivariant}
 \end{matrix}} $.
 \end{itemize}
 Then:
 \begin{enumerate}
 \item \label{AlphaEquiComps-anal}
 $\erg^X$ is analytic,
 and
 \item \label{AlphaEquiComps-equi}
 there is an essentially $(\rho,\tau,\alpha)$-equivariant Borel map
$\phi \colon X \to Y$ if and only if $\erg^X$ is conull
in~$\erg$.
 \end{enumerate}
 \end{prop}

\begin{proof}
  From Prop.~\ref{RohlinDecomp} (and Rem.~\ref{UnionOfProds}), we may
assume the notation of Thm.~\ref{AnalErgComp}.
 Recall that $\rho$, $\tau$, and~$\alpha$ induce Borel maps
 \begin{itemize}
 \item $\rho_{\F} \colon G \times \erg \times \F(S,X) \to \F(S,X)$
\see{BorelIntoF},
 \item $\tau^{\F} \colon \F(S,H) \times \F(S,X) \to \F(S,X)$
\see{PointwiseOps}, and
 \item $\check\alpha \colon G \times \erg \to \F(S,H)$ \see{FiberFunc}.
 \end{itemize}
 Let $G_0$ be a countable, dense subset of~$G$, and define
 $$ \func =  \bigl\{\, (\omega,f ) \in \erg \times \F(S,X)
 \mid \mbox{$\rho_{\F}(g,\omega,f) = \tau^{\F} \bigl(
\check\alpha(g,\omega),f \bigr)$ for all $g \in G_0$} \,\} .$$
 Then $\func$ is an analytic set (in fact, it is Borel, because
$\rho_{\func}$, $\tau^{\func}$, and~$\check\alpha$ are Borel and $G_0$ is
countable).

\pref{AlphaEquiComps-anal} $\erg^X$ is the projection of the analytic
set~$\func$ to~$\erg$.

 (\ref{AlphaEquiComps-equi}~$\Rightarrow$) For a.e.\ $\omega \in \erg$,
the map $\phi_\omega$, defined by $\phi_\omega(s) = \phi(\omega,s)$, is
essentially $(\rho_\omega, \tau, \alpha_\omega)$-equivariant.

(\ref{AlphaEquiComps-equi}~$\Leftarrow$)
 Because $\func$ is analytic, we  may apply the von Neumann Selection
Theorem \pref{vonNeumann}. By assumption, $\erg_{\func}$ is conull
in~$\erg$, so we conclude that there is a Borel function $\phi \colon
\erg \to \F(S,X)$, such that $\bigl( \omega, \phi(\omega) \bigr) \in
\func$, for a.e.\ $\omega \in \erg$. Lemma~\ref{DefdFiberwise} provides
us with a corresponding Borel function $\hat\phi \colon \erg \times S \to
X$. By applying Fubini's Theorem, we see, for each $g \in G_0$, that
 \begin{equation} \label{QuotOfComps-equi}
 \mbox{$\hat\phi \bigl( \rho(g,\omega, s) \bigr) = \tau \bigl(
\alpha(g,\omega,s), \hat\phi(\omega, s) \bigr)$ for a.e.~$(\omega,s) \in
\erg \times S$}.
 \end{equation}
 Let
 $$ H = \{\, g \in G \mid \mbox{\pref{QuotOfComps-equi} holds} \,\} .$$
 Then $H$ is the stabilizer of~$\hat\phi$ under the $\alpha$-twisted
action of~$G$ on $\F(\erg \times S, X)$ (see
Defn.~\fullref{TwistDefn}{action}). Because the action is Borel, we know
that $H$ is a closed subgroup of~$G$. On the other hand, $H$ is dense,
because it contains~$G_0$. Therefore, $H = G$, which means that
$\hat\phi$~is essentially $(\rho,\tau,\alpha)$-equivariant.
 \end{proof}

Any continuous homomorphism $\pi \colon G \to H$ determines a
``constant" cocycle $\pi^\times \colon G \times X \to H$, defined by
$\pi^\times(g,x) = \pi(g)$ \cf{HomoCocDefn}. In the statement of the
following corollary, we ignore the distinction between $\pi$
and~$\pi^\times$.

\begin{cor} \label{QuotErgComps}
  Let
 \begin{itemize}
 \item $G$, $H$, $X$, $\rho$, $\psi$, $\rho_\omega$, $Y$, and $\tau$ be
as in Prop.~\ref{AlphaEquiComps},
 and
 \item $\pi \colon G \to H$ be a continuous homomorphism.
 \end{itemize}
 There exists an essentially $(\rho,\tau,\pi)$-equivariant Borel map $\phi
\colon X \to Y$ if and only if there exists  an essentially
$(\rho_\omega, \tau, \pi)$-equivariant Borel map $\phi_\omega \colon
\psi^{-1}(\omega) \to Y$ for a.e.\ $\omega \in \erg$.
 \end{cor}

\begin{defn}[{\cite[Defn.~9.2.2]{ZimmerBook}}]
 If
 \begin{itemize}
 \item $G$ acts ergodically on $X$,
 \item $\alpha \colon G \times X \to H$ is a Borel cocycle,
 \item $\field$ is a local field,
 and
 \item $H$ is the $\field$-points of an algebraic group over~$\field$,
 \end{itemize}
 then there exists a Zariski-closed subgroup~$L$ of~$H$, such that
 \begin{enumerate}
 \item $\alpha$ is cohomologous to a cocycle taking values in~$L$,
 and
 \item $\alpha$ is \emph{not} cohomologous to a cocycle taking values in
any proper Zariski-closed subgroup of~$H$.
 \end{enumerate}
 The subgroup~$L$ is unique up to conjugacy. It is called
the \emph{algebraic hull} of~$\alpha$.
 \end{defn}

\begin{cor} \label{AlgHull}
 Let
 \begin{itemize}
 \item $(X,\mu)$ be a standard Borel probability space,
 \item $\rho \colon G \times X \to X$ be a Borel action, such that $\mu$
is quasi-invariant,
 \item $\psi \colon X' \to \erg$ be the corresponding ergodic
decomposition \see{ErgDecompThm},
 \item $\field$ be a local field,
 \item $H$ be the $\field$-points of an algebraic group over~$\field$,
 and
 \item $\alpha \colon G \times X \to H$ be a Borel cocycle.
 \end{itemize}
  Then $\alpha$ is cohomologous to a Borel cocyle $\beta \colon G \times
X \to H$, such that, for a.e.\ $\omega \in \erg$, the Zariski closure of
the range of~$\beta_\omega$ is equal to the algebraic hull
of~$\beta_\omega$ {\rm(}where $\beta_\omega$ is the restriction
of~$\beta$ to $G \times \psi^{-1}(\omega)${\rm)}.
 \end{cor}

\begin{proof}
 Choose $\erg'$ as in \fullref{CohoErg}{coc} (with $\beta = \alpha$).
 Chevalley's Theorem \cite[Thm.~5.1, p.~89]{Borel-LinAlgGrps} implies
there is a countable collection $\{(\tau_i,V_i)\}_{i=0}^\infty$ of
rational representations of~$H$, such that every Zariski-closed subgroup
of~$H$ is the stabilizer of some point in some projective space
$\proj{V_i}$. For $c,d, i\in \natural$, define
 \begin{itemize}
 \item $\erg_{c,d} = \bigset{ \omega \in \erg' }{
 \begin{matrix}
 \mbox{the algebraic hull of $\alpha_\omega$ is $d$-dimensional} \\
 \mbox{and has exactly $c$ connected components}
 \end{matrix} }$,
 \item $ Y_{c,d}^i = \bigset{ v \in \proj{V_i} }{
 \begin{matrix}
  \mbox{$\dim \Stab_H(v) = d$, and} \\
 \mbox{$\Stab_H(v)$ has exactly $c$ connected components} \\
 \end{matrix}
 }$,
 \item $Y_{c,d}
= \coprod_{i=0}^\infty Y_{c,d}^i$ (disjoint union),
 and
 \item $\tau \colon H \times Y_{c,d} \to Y_{c,d}$ by $\tau(h,x) = \tau_i(h)
x$ if $x \in \proj{V_i}$.
 \end{itemize}
 Each $Y_{c,d}^i$ is Borel (see \ref{DimBorel} below), so $Y_{c,d}$ is a
standard Borel space.   Therefore, combining the Cocycle Reduction Lemma
\cite[Lem.~5.2.11, p.~108]{ZimmerBook} with
Prop.~\fullref{AlphaEquiComps}{anal} and Rem.~\ref{Anal->Meas} implies
that $\erg_{c,d}$ is absolutely measurable. Thus, we may let
$\erg'_{c,d}$ be a conull, Borel subset of~$\erg_{c,d}$, for each~$d$.

There is no harm in assuming $\erg = \erg'_{c,d}$, for some~$c$ and~$d$.
 Then there is an essentially $(\rho_\omega, \tau,
\alpha_\omega)$-equivariant Borel map from $\psi^{-1}(\omega)$ to
$Y_{c,d}$, for a.e.\ $\omega \in \erg$, so
Cor.~\fullref{AlphaEquiComps}{equi} implies that there is an essentially
$(\rho, \tau, \alpha)$-equivariant Borel map from $X$ to~$Y_{c,d}$. This
means that $\alpha$ is cohomologous to a cocycle~$\beta$, such that the
essential range of~$\beta_\omega$ is contained in an algebraic group
whose dimension and number of connected components are no more than those
of the algebraic hull of~$\beta_\omega$. By changing $\beta$ on a set of
measure~$0$, we may assume that the entire range of~$\beta_\omega$ is
contained in this subgroup. So the desired conclusion follows from the
minimality (and uniqueness) of the algebraic hull of~$\beta_\omega$.
 \end{proof}

The following observation, used in the proof of Cor.~\ref{AlgHull} above,
must be well known, but the authors do not know of a reference.

\begin{lem} \label{DimBorel}
 Suppose
 \begin{itemize}
 \item $G$ is a real Lie group,
 \item $M$ is a Polish space,
 and
 \item $\rho \colon G \times M \to M$ is a continuous action, such that
the stabilizer $\Stab_G(m)$ has only finitely many connected components,
for each $m \in M$.
 \end{itemize}
 Then:
 \begin{enumerate}
 \item \label{DimBorel-dim}
 $\dim \Stab_G(m)$ is a Borel function of $m \in M$.
 \item \label{DimBorel-cpct}
 For each compact subgroup~$K$ of~$G$,
 $$ \{\, m \in M \mid \mbox{$K$ is conjugate to a maximal compact
subgroup of $\Stab_G(m)$} \,\}$$
 is Borel.
 \item \label{DimBorel-dimcpct}
 The dimension of the maximal compact subgroup of $\Stab_G(m)$ is a
Borel function of~$m$.
 \item \label{DimBorel-comp}
 The number of connected components of $\Stab_G(m)$ is a Borel
function of~$m$.
 \item \label{DimBorel-dimcomp}
 For each $d$ and~$c$,
 $$ \bigset{ m \in M }{
 \begin{matrix}
 \mbox{$\dim \Stab_G(m) = d$, and} \\
 \mbox{$\Stab_G(m)$ has exactly~$c$ connected components}
 \end{matrix} }$$
 is Borel.
 \end{enumerate}
 \end{lem}

\begin{proof}
 \pref{DimBorel-dim}
 It is well known (and easy to see) that $\dim \Stab_G(m)$ is an upper
semicontinuous function of $m \in M$.

\pref{DimBorel-cpct} The fixed-point set $M^K$ is closed, and $G$ is
$\sigma$-compact, so $\rho(G, M^K)$ is a countable union of closed sets.
Therefore $\rho(G, M^K)$ is Borel.  Now
 $$ \rho(G, M^K) = \{\, m \in M \mid \mbox{$\Stab_G(m)$
contains a conjugate of~$K$}\,\} ,$$
 so
 $$\bigset{ m \in M }{  \begin{matrix}
 \mbox{$K$ is conjugate to a maximal} \\
 \mbox{compact subgroup of $\Stab_G(m)$}
 \end{matrix} }
 = \rho(G, M^K) \smallsetminus \bigcup_{K' \supset K} \rho(G, M^{K'}) .$$
 Any real Lie group has only countably many conjugacy classes of
compact subgroups \cite[Prop.~10.12]{Adams}, so the union is countable.
Therefore, this is a Borel set.

(\ref{DimBorel-dimcpct}, \ref{DimBorel-comp}) Immediate from
\pref{DimBorel-cpct}. (Recall that the number of components of
$\Stab_G(m)$ is the same as the number of components of any of its
maximal compact subgroups.)

\pref{DimBorel-dimcomp} Combine \pref{DimBorel-dim} and
\pref{DimBorel-comp}.
 \end{proof}


\begin{thebibliography}{GGM}

\bibitem[Ad]{Adams}
 S.~Adams,
 \emph{Reduction of cocycles with hyperbolic targets,}
 Ergodic Theory Dynam. Systems {\bf 16} (1996), no. 6, 1111--1145,
 MR1424391 (98i:58135),
 Zbl 0869.58031.

\bibitem[Ar]{Arveson}
 W.~Arveson,
 \emph{An invitation to $C\sp*$-algebras,}
 Springer-Verlag, New York, 1976,
 MR0512360 (58 \#23621),
 Zbl 0344.46123.

\bibitem[Bo]{Borel-LinAlgGrps}
 A.~Borel,
 \emph{Linear Algebraic Groups,} 2nd ed.,
 Springer, New York, 1991,
 MR1102012 (92d:20001),
 Zbl 0726.20030.

\bibitem[FM]{FisherMargulis}
 D.~Fisher and G.~A.~Margulis,
 \emph{Local rigidity for cocycles,}
 in S.T.Yau, ed.,
 \emph{Surveys in Differential Geometry,}
 Vol. 8,
 International Press, Boston, MA,
 (to appear).

\bibitem[GGM]{GoodrichEtAl}
 K.~Goodrich, K.~Gustafson, and B.~Misra,
	\emph{On converse to Koopman's lemma,}
	 Phys. A {\bf 102}  (1980), no. 2, 379--388,
 MR0582372 (82i:82007).

\bibitem[GS]{Greschonig-Schmidt}
 G.~Greschonig and K.~Schmidt,
 \emph{Ergodic decomposition of
quasi-invariant probability measures,}
 Colloq. Math. {\bf 84/85} (2000), part 2, 495--514,
 MR1784210 (2001i:28021),
 Zbl 0972.37003.

\bibitem[Mc]{Mackey-Ind}
 G.~Mackey,
 \emph{Induced representation of locally compact groups~I,}
 Ann. Math. (2) {\bf 55} (1952) 101--139,
 MR0044536 (13,434a),
 Zbl 0046.11601.

\bibitem[Mr]{MargulisBook}
 G.~A.~Margulis,
 \emph{Discrete subgroups of semisimple Lie groups,}
 Springer-Verlag, Berlin, 1991,
 MR1090825 (92h:22021),
 Zbl 0732.22008.

\bibitem[Ra]{Ramsay}
 A.~Ramsay,
 \emph{Virtual groups and group actions,}
 Adv. Math. {\bf 6} (1971) 253--322,
 MR0281876 (43 \#7590),
 Zbl 0216.14902.

\bibitem[Ro]{Rohlin}
 V.~A.~Rohlin,
 \emph{On the fundamental ideas of measure theory,}
 Amer. Math. Soc. Translation {\bf 10} (1962) 1--54,
 MR0047744 (13,924e).

\bibitem[Ru]{Rudin-RCAnal}
 W.~Rudin,
 \emph{Real and Complex Analysis,}
 2nd ed.,
 McGraw Hill, New York, 1974,
 MR0344043 (49 \#8783),
 Zbl 0278.26001.

\bibitem[WZ]{WheedenZygmund}
 R.~Wheeden and A.~Zygmund,
 \emph{Measure and integral,}
 Marcel Dekker, New York, 1977,
 MR0492146 (58 \#11295),
 Zbl 0362.26004.

\bibitem[Z1]{ZimmerExtn}
 R.~J.~Zimmer,
 \emph{Extensions of ergodic group actions,}
 Illinois J. Math. {\bf 20} (1976) no. 3, 373--409,
 MR0409770 (53 \#13522),
 Zbl 0334.28015.

\bibitem[Z2]{ZimmerBook}
 R.~J.~Zimmer,
 \emph{Ergodic theory and semisimple groups,}
 Birkh\"auser, Basel, 1984,
 MR0776417 (86j:22014),
 Zbl 0571.58015.


\end{thebibliography}
\end{document}